\documentclass{amsart}
\usepackage{amssymb}
\usepackage{latexsym}
\usepackage[dvips]{epsfig}

\newcommand{\bH}{\mathbb H}
\newcommand{\brho}{\bar{\rho}}
\newcommand{\brhoz}{\bar{\rho}_0}
\newcommand{\C}{\mathbb C}
\newcommand{\CP}{\mathbb{CP}^1}

\newcommand{\dpqr}{\Delta(p,q,r)}
\newcommand{\dpqrt}{\Delta(p,q,m)}

\newcommand{\E}{\mathbb{E}}
\newcommand{\fl}[1]{\lfloor #1 \rfloor}
\newcommand{\fltwo}[1]{\lfloor \frac{#1}{2} \rfloor}

\newcommand{\la}{\lambda}
\newcommand{\lt}{\tilde{\lambda}}

\newcommand{\pd}{\pi_1(\St)}
 
\newcommand{\pM}{\partial M}
\newcommand{\pn}{\pi}
\newcommand{\PSLC}{\mbox{PSL}_2(\C)}
\newcommand{\Q}{\mathbb Q}
\newcommand{\R}{\mathbb R}

\newcommand{\SLC}{\mbox{SL}_2(\C)}
\newcommand{\St}{\Sigma_2}
\newcommand{\tXi}{\widetilde{X}_i}

\newcommand{\Z}{\mathbb Z}

\newtheorem{theorem}{Theorem}[section]

\newtheorem{cor}[theorem]{Corollary}
\newtheorem{lemma}[theorem]{Lemma}
\newtheorem{prop}[theorem]{Proposition}

\newcommand{\Pf}{\noindent {\bf Proof: }}
\newcommand{\Rmk}{\noindent {\bf Remark: }}
\newcommand{\Prob}{\noindent {\bf Problem: }}

\begin{document}
\title{Cyclic and finite surgeries on pretzel knots}

\author{Thomas W. Mattman}
\address{Department of Mathematics and Statistics\\ California State
University, Chico\\ Chico CA95929-0525, U.S.A.}
\email{TMattman@CSUChico.edu}
\thanks{Research supported by grants NSERC OGP 0009446 and FCAR EQ 3518.}

\begin{abstract}
We classify Dehn surgeries on $(p,q,r)$ pretzel knots
resulting in a manifold $M(\alpha)$ having cyclic
fundamental group and analyze those leading to a 
finite fundamental group. The proof uses the theory of
cyclic and finite surgeries developed by Culler, Shalen,
Boyer, and Zhang. In particular, Culler-Shalen seminorms 
play a central role.
\end{abstract}

\maketitle

\section{Introduction}
Thurston~\cite{TH} has shown that all but a finite number of 
Dehn surgeries on a hyperbolic knot result in manifolds which 
are again hyperbolic. Those which produce 
a manifold having cyclic or finite fundamental group
are important examples of exceptional (i.e., non-hyperbolic) surgeries.
The theory of Culler-Shalen seminorms has proven to be
a useful tool for understanding these kinds of exceptional surgeries.
These seminorms were first introduced as part of the proof of the Cyclic
Surgery Theorem~\cite{CGLS} and later extended by Boyer and
Zhang~\cite{BZ} to  the study of finite surgeries, eventually leading to
a proof of the Finite Filling Conjecture~\cite{BZ4}.

While those results establish
global bounds on the number of cyclic or finite surgeries
a knot may have, the current paper shows how they
may be refined by focusing on a particular family of knots,
the $(p,q,r)$ pretzel knots. To be specific, we have the 
following theorems.

\begin{theorem} \label{thcyc} The only non-torus pretzel knot 
which admits a non-trivial cyclic surgery is the $(-2,3,7)$
pretzel knot. The non-trivial cyclic surgeries on this knot
are of slope $18$ and $19$.
\end{theorem}

\begin{theorem} \label{thfin}
If a non-torus pretzel knot $K$ admits a non-trivial 
finite surgery, then one of the following holds.
\begin{itemize}
\item $K$ is a $(-2,p,q)$ pretzel knot with $5 \leq p \leq q$ odd
and the filling is not cyclic.
\item $K$ is the $(-2,3,7)$ pretzel knot and the filling is 
  along slope $17,18$, or $19$.
\item $K$ is the $(-2,3,9)$ pretzel knot and the filling is 
  along slope $22$ or $23$.
\end{itemize}
\end{theorem}

Note that a torus knot admits an infinite number of finite
cyclic fillings and that the torus pretzel knots are
well understood~\cite{Kw}.
In particular, if $|p|, |q|,|r| > 1$, then a 
$(p,q,r)$ pretzel knot is torus only if 
$\{p,q,r\} = \{-2,3,3 \}$ or $\{-2,3,5 \}$.
(We consider the $(p,q,r)$ and $(-p, -q, -r)$ knots equivalent as they
are mirror images of one another.)

The surgeries listed for the $(-2,3,7)$ and
$(-2,3,9)$ pretzel knots were discovered by 
Fintushel and Stern~\cite{FS} and
Bleiler and Hodgson~\cite{BH}. The content here
is that there are no other cyclic surgeries
and that the possibilities for a finite 
surgery are restricted. On the other hand,  we know of no
instances of a finite surgery on a knot $(-2,p,q)$ with
$5 \leq p \leq q$. Indeed, we expect that there are none.

Our results can be seen as complementing the 
work of Delman~\cite{D} and others who have been using 
laminations to study Dehn surgery.
In particular, if a hyperbolic arborescent
knot admits a non-trivial cyclic or finite 
surgery, then it must either be a $(p,q,r)$ pretzel
knot or else belong to a certain family of $3$-tangle 
Montesinos knots (see \cite{W2}, especially Theorem
2.1). The current article deals with the first case and
naturally suggests the following:

\vspace{12pt}

\Prob Complete the analysis of cyclic and finite
surgeries on hyperbolic arborescent knots by investigating
Montesinos knots of the from $M(x, 1/p, 1/q)$ where
$x \in \{ -1 \pm 1/2n, -2 + 1/2n\}$ and $p$, $q$, and $n$ 
are positive integers.

\vspace{12pt}

Our main theorems are consequences of the following two results
and work of Delman~\cite{D}.

\begin{theorem} \label{prcyc}
Suppose $K$, a $(-2,p,q)$ pretzel knot ($p, q$ odd and positive),
admits a non-trivial cyclic surgery. Then one
of the following holds. 
\begin{enumerate}
\item $K$ is a torus knot and therefore admits an infinite number of
cyclic surgeries.
In this case either
$\{ p,q \} = \{3,3 \}$ or $\{ p,q \} = \{3,5 \}$ or $\{ p,q \} = \{ 1,n
\}$ for some $n > 0$. 
\item $K$ is the $(-2,3,7)$ pretzel knot and the surgery 
is $18$ or $19$.
\end{enumerate}
\end{theorem}

\begin{theorem} \label{thfinpqr}
 A $(p,q,-r)$ pretzel knot,
with $4 \leq r$ even and $3 \leq p \leq q$ odd admits no
non-trivial finite surgeries.
\end{theorem}

\Pf (of Theorem~\ref{thcyc}) Delman~\cite{D} has shown that if such a knot
admits  a cyclic filling, then it is of the form
$(p,q,-r)$, with $2 \leq r$ even and $3 \leq p \leq q$ odd. 
As only the trivial knot admits a $\Z$ filling \cite{Gb},
Theorem~\ref{thfinpqr} implies further that $r$ must be $2$.
Theorem~\ref{prcyc} completes the proof. \qed

\vspace{12pt}

\Pf (of Theorem~\ref{thfin}) Again, Delman~\cite{D} allows us to reduce
to the case of a $(p,q,-r)$ pretzel knot and Theorem~\ref{thfinpqr}
further shows that $r = 2$. The finite surgeries on $(-2,3,n)$ pretzel
knots are classified in \cite{M1}. That a non-trivial finite
filling of $(-2,p,q)$ with $p \geq 5$ is not cyclic follows from
Theorem~\ref{thcyc}. \qed

\vspace{12pt}

Thus, the main part of this paper is given over to proving
Theorems~\ref{prcyc} and \ref{thfinpqr}. The latter in turn
depends largely on

\begin{theorem} \label{thpqr}
If $K$ is a $(p,q,r)$ pretzel knot with
$p, q$ odd, $r$ even and $1/|p| + 1/|q| + 2/|r| < 1$,
then $K$ admits at most one non-trivial
finite surgery. Moreover such a surgery slope $u$ is odd integral
and there is a non-integral boundary slope in $(u-1,u+1)$. 
\end{theorem}

Theorems~\ref{prcyc}, \ref{thpqr}, and \ref{thfinpqr}
will be proved in Sections~\ref{seccyc}, \ref{secone}, and 
\ref{secfin}. Section~\ref{secCS}, which follows, introduces notation
and provides a brief review of Culler-Shalen theory
which will play a central role in our arguments.

\section{\label{secCS}%
Notation and Culler-Shalen Theory}

\begin{figure}
\begin{center}
\epsfig{file=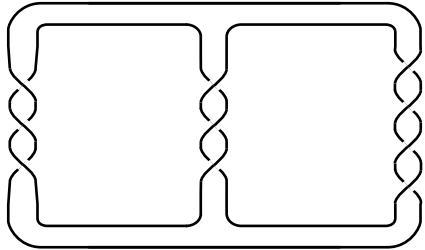}
\caption{\label{fg334}%
The $(-3,3,4)$ pretzel knot.}
\end{center}
\end{figure}

In this section let $K$ denote a $(p,q,r)$ pretzel knot.
Our sign conventions are illustrated by Figure~\ref{fg334}
which shows the $(-3,3,4)$ pretzel knot. Pretzel knots
whose indices $p,q,r$ agree up to permutation are
ambient isotopic. Moreover, taking the mirror reflection
corresponds to changing the signs of all the indices. As
this reduces to an isomorphism of the knot group $\pi$,
we will consider the knots $(p,q,r)$ and
$(-p,-q,-r)$ equivalent.

By \cite{O}, $K$ is 
{\em small} in the sense that the knot complement $M = S^3 \setminus K$
contains no closed essential surfaces. The knot is therefore
either torus or hyperbolic. The torus pretzel knots are
classified in \cite{Kw}. 

The fundamental group of the $2$-fold branched cyclic cover
$\pd$ of $K$ is a central extension of the triangle group $\dpqr$.
Let $\mu$ denote the class of a meridian in $\pi$ and
$\la$ that of a preferred longitude.
We will use $\{ \mu, \la \}$ coordinates to identify the 
surgery slopes on $K$ with $\Q \cup \{1/0 \}$ and denote
Dehn surgery along slope $a/b$ by $M(a/b)$.
Surgery along the meridian $M(1/0) = M(\mu) = S^3$ is called
{\em trivial} surgery.
The {\em distance} between two surgery slopes $a/b$ and $c/d$ 
is the minimal geometric intersection number 
$\Delta(a/b, c/d) = |ad - bc|$.
We will say that $a/b$ is a {\em boundary slope} if there
is an essential surface in $M$ which meets $\pM$ in a non-empty
set of curves having slope $a/b$. An {\em essential surface} is
one which is properly embedded, orientable, incompressible,
$\partial$-incompressible, and non-$\partial$-parallel.

We now briefly introduce Culler-Shalen theory under
the assumption that $K$ is a small, hyperbolic knot.
A more detailed account may be found in 
\cite[Chapter 1]{CGLS} and \cite{BZ2}.

Let $R = \mbox{Hom}(\pn, \SLC)$ denote the set of
$\SLC$-representations of the fundamental group of $M$.
Then $R$ is an affine algebraic set, as is $X$, the set of characters
of representations in $R$.
Since $M$ is small, the irreducible components of $X$ are 
curves~\cite[Proposition 2.4]{CCGLS}. Moreover, for each component
$R_i$ of $R$ which contains an irreducible representation, the 
corresponding curve $X_i$ induces a non-zero seminorm $\| \cdot \|_i$ on 
$V = H_1(\pM; \R)$~\cite[Propositon 5.7]{BZ2} via the following construction.

For $\gamma \in \pn$, define the regular function
$I_{\gamma}:X \to \C$ by $I_{\gamma}(\chi_{\rho}) = \chi_{\rho}(\gamma) =
\mbox{trace}(\rho(\gamma))$.
By the Hurewicz isomorphism,
a class $\gamma \in L =  H_1(\pM ;{\Z})$ determines an element
of $\pi_1(\pM)$, and therefore an element of $\pn$ well-defined
up to conjugacy. The function
$f_{\gamma} = I_{\gamma}^2 -4$ is again regular and so can be pulled back to
$\tXi$, the smooth projective variety birationally equivalent to
$X_i$. For $\gamma \in L$, $\| \gamma \|_i$ is the degree of
$f_{\gamma} : \tXi \to \CP$. In practice, this degree can
often be calculated using $Z_x(f_{\gamma})$, the degree
of zero of $f_{\gamma}$ at a point $x$ in the character variety.
 The seminorm is extended to $V$ 
by linearity. We will call a seminorm constructed in this manner 
a Culler-Shalen seminorm.

Let $\| \cdot \|_T$ denote the sum of the 
Culler-Shalen seminorms, i.e.,  $\| v \|_T = \sum_i \|
v  \|_i$ (here $v \in V$) and let 
$S = \min \{ \| \gamma \|_T \, ; \, \gamma \in L, \, \| \gamma \|_T > 0
\}$  be the minimal norm. 
Note that, as $K$ is hyperbolic,
$\| \cdot \|_T$ will be a norm (and not just a seminorm) 
\cite[Chapter 1]{CGLS}.

\section{\label{seccyc}%
Proof of Theorem~\ref{prcyc}}

In this section, let $K$ be a $(-2,p,q)$ pretzel knot 
($p,q$ odd and positive, $p \leq q$).
The following lemma can be proved using the methods
of Hatcher and Oertel~\cite{HO}.

\begin{lemma} \label{le2pq}
Assume $3 \leq p \leq q$.
If $p \geq 7$ (respectively $q \geq 7$) then 
$$\frac{p^2 - p - 5}{\frac{p-3}{2}} \quad \quad \quad (\mbox{resp. }
\frac{q^2 - q - 5}{\frac{q-3}{2}})$$
is a non-integral boundary slope of $K$. Moreover, these are 
the only non-integral boundary slopes of $K$.
\end{lemma}

\begin{lemma} \label{letor}
$M(2(p+q))$ contains an incompressible torus. \end{lemma}

\Pf
Use the double cover of the ``obvious''
spanning surface of the knot (which is a punctured Klein bottle). \qed

\setcounter{section}{1}
\setcounter{theorem}{2}

\begin{theorem} 
Suppose $K$
admits a non-trivial cyclic surgery. Then one
of the following holds. 
\begin{enumerate}
\item $K$ is a torus knot and therefore admits an infinite number of
cyclic surgeries.
In this case either
$\{ p,q \} = \{3,3 \}$ or $\{ p,q \} = \{3,5 \}$ or $\{ p,q \} = \{ 1,n
\}$ for some $n > 0$. 
\item $K$ is the $(-2,3,7)$ pretzel knot and the surgery 
is $18$ or $19$.
\end{enumerate}
\end{theorem}

\setcounter{section}{3}
\setcounter{theorem}{2}

\Pf
Theorem III of \cite{Kw} shows that $K$ is torus iff it is 
as characterized in 1.

Since $K$ is small \cite{O}, we can assume that
$K$ is hyperbolic. The cyclic surgeries $18$ and $19$ of the 
$(-2,3,7)$ pretzel knot were first observed by Fintushel and 
Stern (see \cite[Section 4]{FS}). 
Our task is to show that there is no other choice for
$p$ and $q$ leading to a cyclic surgery.

The case $p=3$ is the subject of \cite{M1} where we  
show that there are no non-trivial cyclic surgeries
when $q \geq 9$ and that the cyclic surgeries of the
$(-2,3,7)$ pretzel knot are as stated. 

If $p=5$, the boundary slopes \cite{HO} are $0,14,15, 
\frac{q^2-q-5}{\frac{q-3}{2}}, 2q+10$, and $2q+12$. By
\cite[Theorem 4.1]{Du}, a non-trivial cyclic surgery could occur
only at $2q+4$ or $2q+5$. However, as we explain below, a
cyclic surgery would have to be within distance $5$ of the 
toroidal surgery $2q+10$ (see Lemma~\ref{letor}). So the only
candidate is $2q+5$. Now the $(-2,5,5)$ pretzel has no non-integral
boundary slopes so (\cite[Theorem 4.1]{Du}) it has no non-trivial
cyclic surgeries. As for $(-2,5,7)$, SnapPea~\cite{Wk} 
shows that $2q+5 = 19$ surgery on this knot is hyperbolic. So
we can assume $q \geq 9$.

Suppose (for a contradiction) that $2q+5$ is indeed a cyclic
surgery. By ~\cite[Lemma 6.2]{BZ}, the (total) norm can be written
\begin{eqnarray*} 
\| \gamma \|_T & = & 2 {[} a_1 \Delta( \gamma, 0) + 
a_2 \Delta( \gamma, 14) + a_3 \Delta( \gamma, 15) \\
&&\mbox{ } + a_4 
\Delta( \gamma, \frac{q^2-q-5}{\frac{q-3}{2}}) +
a_5 \Delta( \gamma, 2q + 10) + a_6 \Delta( \gamma, 2q + 12) {]}
\end{eqnarray*}
where the $a_i$ are non-negative integers.
If $2q+5$ is cyclic it has minimal norm $S$, as does the meridian
surgery $\mu$ (\cite[Corollary 1.1.4]{CGLS}).
The norm of $2q+4$ will also be of
interest, and it will be bounded by the minimal norm $S$.
\begin{eqnarray*}
S & = & \| \mu \|_T = 2 {[}a_1 + a_2 + a_3 + \frac{q-3}{2} a_4 + a_5 +
a_6 {]} \\
S & = & \| 2q + 5 \|_T = 2 {[} (2q + 5) a_1 + (2q-9) a_2 + (2q-10) a_3
+ \frac{q-5}{2} a_4 + 5 a_5 + 7 a_6 {]} \\
S & \leq & \| 2q+4 \|_T =
2{[} 2q+4 a_1 + (2q-10) a_2 + (2q-11) a_3 + a_4 + 
+ 6 a_5 + 8 a_6 {]} 
\end{eqnarray*}

Subtracting the first two equations, we have
\begin{equation}
a_4 = (2q+4)a_1 + (2q-10) a_2 + (2q-11) a_3 + 4a_5 + 6 a_6,
\label{eqa4}
\end{equation}
while subtracting the second from the third leaves
\begin{eqnarray*}
&&a_5 + a_6 \geq a_1 + a_2 + a_3 + \frac{q-7}{2} a_4, \\
& \Rightarrow & \eta (a_5 + a_6 - a_1 - a_2 -a_3) \geq a_4,
\end{eqnarray*}
where $\eta = \frac{2}{q-7} \leq 1$.

Combining this with Equation~\ref{eqa4}, we have
$$ 0 \geq (2q+4 + \eta) a_1 + (2q-10 + \eta) a_2 + (2q-11 +
\eta)a_3 + (4-\eta) a_5 + (6-\eta) a_6.$$
Since $a_i \geq 0$, this shows $a_1 = a_2 = a_3 = a_5 = a_6 = 0$. 
On the other hand, for a norm, at least two of the $a_i$ must
be non-zero. This contradiction shows that there can be no non-trivial 
cyclic surgery when $p=5$.

So let us assume $7 \leq p \leq q$. Dunfield\cite[Theorem 4.1]{Du}
has shown that any non-trivial cyclic surgery on a knot such as $K$ must
lie near a non-integral surgery. Combining this with
Lemma~\ref{le2pq}, the only candidates for a non-trivial cyclic 
surgery are $2p +4$, $2p+5$, $2q + 4$, and $2q+5$. Suppose that 
$u$ is one of these candidates slopes and $M(u)$ is a cyclic filling.
Since $K$ is strongly invertible, 
the Orbifold Theorem implies that $M(u)$ admits
a geometric decomposition (see \cite[Corollary 1.21]{CHK}).
Now, as $\Delta(u,2(p+q)) > 5$, $M(u)$ is irreducible \cite{Oh,Wu} and 
atoroidal \cite{G1} and therefore has a geometric structure.

Note that $\pi_1(M(u)) \ncong \Z$ (see \cite{Gb}), so 
$\pi_1(M(u))$ is finite. The geometry is therefore $S^3$, and as 
$\pi_1(M(u))$ is finite cyclic, we deduce that $M(u)$ is a lens
space. However, this contradicts \cite[Theorem 1.1]{G2} which 
states that the distance between a lens space surgery such as $u$
and a toroidal surgery such as $2(p+q)$ is at most $5$.
We conclude that there are also no non-trivial cyclic surgeries in 
this case. \qed

\section{\label{secone}%
Proof of Theorem~\ref{thpqr}}

Let $K$ be a $(p,q,r)$ pretzel knot where
$p = 2k+1$, $q = 2l+1$ and $r = 2m$. We will be assuming
that $1/|p| + 1/|q| + 1/|m| \leq 1$, and
this ensures that $K$ is hyperbolic~\cite{Kw}.  

\begin{lemma} \label {leinfev}
Let $1/|p| + 1/|q| + 1/|m| \leq 1$.
If $b$ is odd, then the Dehn filling
$M(2a/b)$ of the knot complement $M$ has infinite $\pi_1$.
\end{lemma}

\Pf
As $\dpqrt$ is infinite, our strategy is 
to construct a representation of $\pi_1(M(2a/b))$ with image
$\dpqrt$. Changing the sign of $p$, $q$, or $r$ will
not change the triangle groups $\dpqrt$ and $\dpqr$, so we will assume 
$p,q,r > 0$ in order to simplify the notation. For the
general case, use $|p|$, $|q|$, and $|r|$ instead.

First note that there is a faithful $\PSLC$-representation 
of $\dpqrt$. Indeed, either $\{p,q,m \} = \{3,3,3 \}$, and 
$\dpqrt$ is a set of isometries of the Euclidean
plane $\E^2$, or else
$\frac{1}{p} + \frac{1}{q} + \frac{1}{m} < 1$
and $\dpqrt$ represents isometries of the hyperbolic plane $\bH^2$.
Now, both $\E^2$ and $\bH^2$ imbed isometrically into hyperbolic
$3$-space $\bH^3$. For example, in the upper half space model,
the set $\{ z = 1 \}$ is a Euclidean plane, while the $xz$ plane is 
$\bH^2$. Moreover, isometries of these planes are restrictions
of isometries of $\bH^3$. Thus, in either case, $\dpqrt$ embeds in 
$\PSLC$, the set of orientation preserving 
isometries of $\bH^3$. This provides the required
faithful representation.

Let $\brhoz$ be a representation of $\dpqr$ obtained by composing
the obvious homomorphism $\dpqr \to \dpqrt$ with a faithful
representation of $\dpqrt$ in $\PSLC$.
Then (as in \cite[Proposition 1.1]{M1}) 
$\brhoz$ ``extends'' to a $\PSLC$-representation $\brho$
of the knot group $\pi$ which in turn lifts to an $\SLC$-representation
$\rho$. (The obstruction to such a lift is in $H^2(\pi; \Z / 2)$ 
\cite[Section 3]{BZ2}. For a knot in $S^3$, the second cohomology
is trivial and there is no obstruction.) 
Moreover, $\brho(\mu^2) = 1$.

On the other hand, we can determine
the image of $\lambda$ in $\dpqr =  \langle f,g,h \mid f^{r}, g^{p},
h^{q}, fgh \rangle$ to be $\lt = g^k f^m g^{k+1} h^l f^m h^{l+1}$
(compare \cite{Tr}).
Now, as $\brhoz$ factors through a representation of $\dpqrt$, 
we have $\brhoz(f^m) = 1$ and consequently $\brhoz(\lt) = 1$.
Then $\brho(\lambda) = 1$ as well. 

So, for any filling of the form $\alpha = 2a/b$, we have 
$\brho(\alpha) = 1$ whence  
$\brho$ factors through $\pi_1(M(\alpha))$. 
Since $\brho$ is an extension of $\brhoz$, which factors through a
faithful representation  of the infinite group $\dpqrt$,
we see that $\pi_1(M(\alpha))$ must also
be  infinite. \qed

So under the hypothesis of the lemma, every
$2a/b$ filling of $K$ is infinite. This means
that any finite surgeries would have to be of the form $(2a+1)/b$ 
and therefore would have norm $\|(2a+1)/b\|_T \leq S + 8$ 
\cite[Theorem 2.3]{BZ}. On the other hand, the $2a/b$ fillings
will have norm larger than $S+8$. 

\begin{lemma} \label{lebigev}  Let
$1/|p|+ 1/|q|+1/|m| < 1$. If $b$ is odd, then $\|2a/b \|_{T} \geq S +
12$.
\end{lemma}

\Pf Again, we assume $p,q,r > 0$.
 
We first observe that there are at least three irreducible
$\PSLC$-characters of $\dpqrt$. Indeed, 
by \cite[Proposition D]{BB},
the number of $\PSLC$-characters of $\Delta(p,q,r)$ is
\begin{eqnarray} \label{eqpqr}
&& (p- \fltwo{p} - 1)(q- \fltwo{q} - 1)(r- \fltwo{r} - 1) +
   \fltwo{p} \fltwo{q} \fltwo{r} \\
&& \mbox{}+ \fltwo{\gcd(p,q)} + \fltwo{\gcd(p,r)} + \fltwo{\gcd(q,r)}
   + 1 \nonumber \end{eqnarray}
where $\fl{x}$ denotes the largest integer less than or equal to $x$.
This count includes characters of reducible representations.
The character of 
a reducible representation is also the character of 
a diagonal (hence abelian) representation. 
So, to count the characters of reducible representations 
we can look at representations of $H_1(\Delta(p,q,r)).$
Let $a = \gcd(p,q,r),$ $b = \gcd(pq,pr,qr).$ Then
$H_1(\Delta(p,q,r)) = \Z/a \oplus \Z/(b/a)$
and hence $|H_1(\Delta(p,q,r))| =
b.$ Consequently, the number of characters of $H_1(\Delta(p,q,r))$ is
\begin{equation} \label{eqH1}
   \begin{array}{cl} \fltwo{b} + 1, & \mbox{if } a \equiv 1 \pmod{2},  \\
                     \fltwo{b} + 2, & \mbox{if } a \equiv 0 \pmod{2}.
   \end{array}
\end{equation}

Thus, by taking the difference of (\ref{eqpqr}) and
(\ref{eqH1}), we see
that $\dpqrt$ has at least three irreducible
$\PSLC$-characters. Using the method of the previous lemma,
these can be used to construct irreducible $\PSLC$-characters
of $\pi_1(M(\alpha))$ when $\alpha$
is of the from $2a/b$. None of these characters
are dihedral, so each is covered twice in 
$\SLC$ (see \cite[Lemma 5.5]{BZ}). As they are the characters of
irreducible  representations of a triangle group, 
they are smooth points of $X(M)$ (see \cite[Proposition 7]{BZ5}).
 
Moreover, they are zeroes of $f_{\alpha}$ which are not zeroes
of $f_{\mu}$. (As in the previous lemma, these are characters of
representations which take $\mu$ to an element of order two.)
It follows from \cite[Theorem A]{BB} that $Z_x(f_{\alpha}) 
= Z_x(f_{\mu}) + 2$, and since we have six such characters $x$, we see 
that $\| 2a/b \|_T = \| \alpha \|_T \geq \| \mu \|_T + 12 = S + 12$.
\qed

\setcounter{section}{1}
\setcounter{theorem}{4}

\begin{theorem} 
If $K$ is a $(p,q,r)$ pretzel knot with
$p, q$ odd, $r$ even and $1/|p| + 1/|q| + 2/|r| < 1$,
then $K$ admits at most one non-trivial
finite surgery. Moreover such a surgery slope $u$ is odd integral
and there is a non-integral boundary slope in $(u-1,u+1)$. 
\end{theorem}

\setcounter{section}{4}
\setcounter{theorem}{2}

\Pf 
The conditions on $p,q,r$ ensure that $K$ is hyperbolic 
\cite{Kw}.

Let $\alpha$ be a finite surgery of such a knot. We have 
already observed (Lemma~\ref{leinfev}) that $\alpha = (2a+1)/b$.
Since meridional surgery is cyclic, we can apply \cite[Theorem 1.1]{BZ}
to see that $b \leq 2$. 

If $\alpha = (2a+1)/2$ were a finite filling, then, by \cite[Theorem 2.3]{BZ},
$\| \alpha \|_T \leq S + 8$. At the same time, $\| \mu \|_T = \| - \mu \|_T 
= S$. The line joining $\alpha = (2a+1,2)$ and $\mu = (1,0)$ in surgery space 
$V \cong H_1(\pM ; \R) \cong \R^2$ passes through $(a+1,1)$ while
the line through $\alpha$ and $ - \mu$ passes through $(a,1)$.
It follows that $\|a+1\|_T$ and $\|a\|_T$ are both less than $S+4$.
Since one of them is even, this contradicts Lemma~\ref{lebigev}.

So any non-trivial finite fillings must be odd integral. Suppose there were 
two such. Each would have norm at most $S+8$. The line joining them
would necessarily pass through some even integral surgeries which 
would therefore also have norm at most $S+8$. This again
contradicts Lemma~\ref{lebigev}.

Now suppose that $2a+1$ is a non-trivial finite filling.
Then $\| 2a+1,1 \|_T \leq S+8$ while $\|2a,1\|_T \geq S+12$
by Lemma~\ref{lebigev}. Let $P \subset V$ denote the norm-ball of
radius $S+8$.
By \cite[Proposition 1.1.2]{CGLS}, $P$ is a finite-sided
convex polygon whose vertices are multiples of boundary slopes. 
In particular, $(2a+1,1)$ is not a vertex of $P$ 
(otherwise \cite[Theorem 2.0.3]{CGLS}, $M(2a+1) \cong S^2 \times S^1$
which is absurd).  

We now construct the non-integral boundary slope $c/d$ and show
that it lies in the interval $(2a, 2a+2)$.
Since $(2a+1,1)$ is inside $P$ and $(2a,1)$ is not, there is a segment of
$\partial P$ which intersects the line $y=1$ between them. Let $k(c,d)$ be
the vertex of this segment which lies on or above $y=1$, i.e., $k \in \Q$
and $c/d$ is a boundary slope.  Consider the segment
from the origin to $k(c,d)$. As both endpoints are in $P$, this segment
is also. It crosses $y=1$ at $(c/d,1)$ which must lie between $(2a,1)$
and $(2(a+1),1)$. (Otherwise, the segment joining $(2a+1,1)$ and
$(c/d,1)$ passes through $(2a,1)$, say. Since both endpoints are in 
$P$, this segment is in $P$ and in particular $(2a,1)$ is in $P$, a
contradiction.)

Thus $|2a+1  - c/d| < 1$, as required.  \qed

\begin{cor} \label{cornonint}
If a knot satisfies the conditions of the 
theorem and has no non-integral boundary slopes, then
it admits no non-trivial finite surgeries.
\end{cor}

\begin{cor} Alternating  $(p,q,r)$ pretzel knots with
$1/|p| + 1/|q| + 2/|r| < 1$ admit no non-trivial finite surgeries.
\end{cor}

\Pf This follows since pretzel knots are Montesinos knots and alternating
Montesinos knots have no non-integral boundary slopes (see~\cite{HO}).
\qed

\vspace{12pt} 

\Rmk
Note that the second Corollary also follows from Delman and
Robert's~\cite{DR} proof that alternating knots satisfy strong property
P. 

\section{\label{secfin}%
Proof of Theorem~\ref{thfinpqr}}

We turn now to the case of a $(p, q, -r)$ pretzel knot $K$
where $4 \leq r$ is even and $p$ and $q$ are both
odd. We will assume $3 \leq p \leq q$. The strategy is similar
to that of Section~\ref{seccyc}. We begin with two lemmas based on the
work of Hatcher and Oertel~\cite{HO}.

\begin{lemma} \label{lenoint}
If $p \geq 2r+1$, then
\begin{equation} \label{eqbdy}
\frac{p(p-1)+1 -3r}{\frac{p-1-r}{2}} \quad \quad
\mbox{ and } \quad \quad \frac{q(q-1)+1 -3r}{\frac{q-1-r}{2}}
\end{equation}
are the non-integral boundary slopes of $K$.
\end{lemma}

\begin{lemma} \label{lepsmall}
If $p < r$, and $K$ has a non-integral boundary slope, then 
that slope is
\begin{equation} \label{eqpsmall}
2(p+q+r-1 - \frac{(p-1)(q-1)}{p-1 + q-1}). \end{equation}
\end{lemma}

The proof of the following lemma is similar to that of 
Lemma~3.2.

\begin{lemma} \label{letor2}
$M(2(p+q))$ contains an incompressible torus.
\end{lemma}

\begin{prop} \label{prplarge}
If $p > 2r+1$, then $K$ admits no non-trivial finite
surgeries. 
\end{prop}

\Pf 
By Theorem~\ref{thpqr}, a non-trivial finite surgery would be 
close to one of the non-integral boundary slopes of
Lemma~\ref{lenoint}. However,
\begin{eqnarray*}
2(p+q) - \frac{p(p-1)+1 -3r}{\frac{p-1-r}{2}} & = & 
2(p+q) - 2(p+r) - \frac{(r-1)^2}{\frac{p-1-r}{2}}\\
& = & 2q - 2r - \frac{(r-1)^2}{\frac{p-1-r}{2}} \\
& \geq & 4r+6 - 2r - \frac{(r-1)^2}{\frac{r}{2}+1} \\
& = & \frac{7r+5}{\frac{r}{2}+1} \geq 11 
\end{eqnarray*}
and similarly for the other slope of Lemma~\ref{lenoint}.
Therefore, any non-trivial finite surgery would be of distance
(in the sense of minimal geometric intersection) 
greater than $10$ from the toroidal surgery $2(p+q)$. However
this contradicts work of 
Agol~\cite{Ag} and Lackenby~\cite{L}
showing that the distance between exceptional surgeries
is $\leq 10$. \qed 

\begin{prop} \label{prpsmall}
If $p \leq r-5$, then $K$ admits no non-trivial finite 
surgeries.
\end{prop}

\Pf As in the previous proposition, we observe that 
$$| 2(p+q+r-1 - \frac{(p-1)(q-1)}{p-1 + q-1}) - 2(p+q) | > 10.$$
Thus the lone non-integral boundary slope of Lemma~\ref{lepsmall}
is too far from the toroidal boundary slope $2(p+q)$
(by Theorem~\ref{thpqr} a finite filling could only occur at
an odd-integral slope, which would therefore have to be within
distance $9$ of the even number $2(p+q)$). 
\qed

We now have a fairly precise description of what a finite
filling $s$ on a $(p,q,-r)$ pretzel knot would look like.
By Theorem~\ref{thpqr}, $s$ would have to be odd-integral
and near a non-integral boundary slope and by 
Propositions~\ref{prplarge} and \ref{prpsmall}, 
$p+3 \geq r \geq (p-1)/2$. We now propose to 
explicitly calculate the fundamental group of such a filling.
We will then project onto a group $G$ and observe that 
$G$ is generically infinite. 

The Wirtinger 
presentation~\cite[Section 3.D]{R} of a $(p,q,-r)$ pretzel knot is 
(compare \cite[Equation 1]{Tr}):
\begin{eqnarray*}
\pi_1(M) = \langle x,y,z &\mid&
(zx)^{(p-1)/2}z(zx)^{(1-p)/2} = (yx)^{-(q+1)/2}y(yx)^{(q+1)/2},\\
&&(yz^{-1})^{-r/2}y(yz^{-1})^{r/2} = (yx)^{(1-q)/2}x(yx)^{(q-1)/2},\\
&&(yz^{-1})^{-r/2}z(yz^{-1})^{r/2} = (zx)^{(p+1)/2}x(zx)^{-(p+1)/2}
\rangle .
\end{eqnarray*}
The longitude being 
$$l = x^{-2(p+q)}(yx)^{(q-1)/2}(yz^{-1})^{-r/2}(yx)^{(q+1)/2}(zx)^{(p-1)/2}
(yz^{-1})^{r/2}(zx)^{(p+1)/2},$$
filling along an odd integral slope $s$ results in 
\begin{eqnarray*}
\pi_1(M(s)) = \langle x,y,z &\mid&
(zx)^{(p-1)/2}z(zx)^{(1-p)/2} = (yx)^{-(q+1)/2}y(yx)^{(q+1)/2},\\
&&(yz^{-1})^{-r/2}y(yz^{-1})^{r/2} = (yx)^{(1-q)/2}x(yx)^{(q-1)/2},\\
&&(yz^{-1})^{-r/2}z(yz^{-1})^{r/2} = (zx)^{(p+1)/2}x(zx)^{-(p+1)/2}, x^{s}l
\rangle .
\end{eqnarray*}

We can obtain a more manageable factor group $G$ by adding
the relators $(yz^{-1})^{r/2}$, $yx^{-1}$, and $(zx)^p$:
\begin{eqnarray*}
G &=& \langle y,z \mid (yz^{-1})^{r/2}, (zy)^p, 
z = (zy)^{(p+1)/2}y(zy)^{-(p+1)/2}, y^{s-2p}
\rangle \\
& = & \langle w,y \mid (y^2w^2)^{r/2}, w^p, (wy)^2, y^{s-2p} \rangle,
\end{eqnarray*}
where $w = (zy)^{(p-1)/2}$. This is an example of a group
which Coxeter~\cite{C} has called
$$ ( 2,a,b;c ) = \langle R,S \mid R^a, S^b, (RS)^2,
(R^2S^2)^c  \rangle.$$
Thus $G = ( 2,p,|s-2p|;r/2 )$. Moreover, $\pi_1(M(s))$ will
be infinite whenever $G$ is.

And indeed, these groups are usually infinite as Edjvet has shown:
\begin{theorem}[Main Theorem of \cite{E}] \label{thE}
If $2 \leq a \leq b$, $2 \leq c$ and $(2,a,b;c) \neq (2,3,13;4)$,
then the group $(2,a,b;c)$ is finite if and only if it is one of
the following:
\begin{tabbing}
  (viii) \= $(2,2,b;c)$ \= $(2 \leq b, 2 \leq c)$ \kill
  (i)   \> $(2,2,b;c)$ \> $(2 \leq b, 2 \leq c)$; \\
  (ii)  \> $(2,3,b;c)$ \> $(3 \leq b \leq 6, 4 \leq c)$; \\
  (iii) \> $(2,3,7;c)$ \> $(4 \leq c \leq 8)$; \\
  (iv)  \> $(2,3,b;c)$ \> $(8 \leq b \leq 9, 4 \leq c \leq 5)$; \\
  (v)   \> $(2,3,b;4)$ \> $(10 \leq b \leq 11)$;\\
  (vi)  \> $(2,4,b;2)$ \> $(4 \leq b)$; \\
  (vii) \> $(2,4,4;c)$ \> $(3 \leq c)$; \\
  (viii) \> $(2,4,5;c)$ \> $(3 \leq c \leq 4)$; \\
  (ix)  \> $(2,4,7;3)$; \\
  (x)   \> $(2,5,b;2)$ \> $(5 \leq b \leq 9)$;\\
  (xi) \> $(2,6,7;2)$.
\end{tabbing}
\end{theorem}

Since $p$ is odd, if $p \leq |s-2p|$, then $G$ is infinite unless 
$p = 3$ or $5$.

Similarly, if $|s-2p| \leq p$, we
see that $G$ is infinite unless $|s-2p|= 3$ or $5$ ($|s-2p|$ is also
odd) whence $s \leq 2p + 5$.
 On the other hand, by \cite{Ag,L}, the finite filling $s$  
and the toroidal filling $2(p+q)$ (Lemma~\ref{letor2}) have distance 
at most
$10$. Since $s$ is odd and $2(p+q)$ is  even, they in fact differ by at
most
$9$. Thus,
$9 \geq 2(p+q) - s \geq 2(p+q) - (2p+5) = 2q-5$.
It follows that $2q \leq 14$, whence $3 \leq p \leq q \leq 7$.

So we can assume that $3 \leq p \leq 7$. That is, if $p \geq 9$,
the knot admits no non-trivial finite surgeries. 

Now, earlier work shows that there are no non-trivial surgeries
unless $\frac{p-1}{2} \leq r \leq p+3$. So,   
given $r \geq 4$, and assuming $3 \leq p \leq 7$, 
we see that we are left to investigate $4 \leq r \leq 10$.
And since Theorem~\ref{thpqr} does not apply to the knots $(-4,3,3)$, 
$(-4,3,5)$ and $(-6,3,3)$ these knots must also be 
examined. The details may
be found in \cite{M3}. Note that 
we again use SnapPea~\cite{Wk} to resolve a few difficult cases.

In summary then, Theorems~\ref{thpqr} and \ref{thE}
combine to show that a $(p,q,-r)$ admits no non-trivial finite
surgery unless $4 \leq r \leq 10$ and $3 \leq p \leq 7$. 
Those cases have been investigated directly and 
also admit no non-trivial finite surgeries. We have
therefore proved 

\setcounter{section}{1}
\setcounter{theorem}{3}

\begin{theorem} 
A $(p,q,-r)$ pretzel knot,
with $4 \leq r$ even and $3 \leq p \leq q$ odd admits no
non-trivial finite surgeries.
\end{theorem}

\section*{Acknowledgments}
This work forms part of my Ph.D. thesis and
I would like to thank my supervisor Steven Boyer for
his substantial contributions and indispensable advice.
I am grateful to Nabil Sayari and Xingru Zhang for helpful conversations
and to Andreas Caranti and Rick Thomas
for telling me of Edjvet's paper.

\end{document}